\documentclass[11pt]{amsart}
\usepackage{amsfonts} %za izpis oznak za mnozice stevil
\usepackage{amsmath} %za definicijo matematicnih operatorjev, kot je npr. arctg, in nekaterih simbolov
\usepackage{amsthm} %izpise piko za stevilko definicije, ce je definicija vpeljana kot newtheorem
\usepackage{amssymb} % za razne simbole
\usepackage{enumitem} % da lahko nadaljujes enumerate

\newcommand{\sep}{\itemsep 0pt}
\newcommand{\set}[2]{\{#1 \ ; \ #2\}}
\newcommand{\ideal}{\triangleleft}

\DeclareMathOperator{\rank}{rank}
\DeclareMathOperator{\soc}{soc}
\DeclareMathOperator{\End}{End}

\newtheorem{theorem}{Theorem}[section]
\newtheorem{proposition}[theorem]{Proposition}
\newtheorem{lemma}[theorem]{Lemma}
\newtheorem{corollary}[theorem]{Corollary}

\theoremstyle{definition}
\newtheorem{definition}[theorem]{Definition}
\newtheorem{example}[theorem]{Example}
\newtheorem{remark}[theorem]{Remark}

\begin{document}

\title[Generalized corner rings]{Structure theorem for generalized corner rings}

\author{Nik Stopar}
\address{Faculty of Electrical Engineering, University of Ljubljana, Tr\v za\v ska cesta 25, 1000 Ljubljana, Slovenia}
\email{nik.stopar@fe.uni-lj.si}

\begin{abstract}
We apply recent results on the rank of elements of rings to study the structure of generalized corner rings $aRa$, where $R$ is a unital ring and $a$ an element of $R$. We give a complete description of the structure of $aRa$ when $a^2$ has finite rank and provide an example to show that this assumption is necessary and optimal.
\end{abstract}

\maketitle

{\footnotesize \emph{Key Words:} corner ring, direct sum decomposition, rank, regular element

\emph{2010 Mathematics Subject Classification:}
Primary
%16D25, % Modules, bimodules and ideals - Ideals
16D70, % Modules, bimodules and ideals - Structure and classification, direct sum decomposition, cancellation
16G99, % Representation theory of rings and algebras - None of the above, but in this section
16U99. % Conditions on elements - None of the above, but in this section
%Secondary
%16N60, % Radicals and radical properties of rings - Prime and semiprime rings
%16E50. % Homological methods - von Neumann regular rings and generalizations
}

\section{Introduction}

The aim of this paper is to present a particular application of the results on the rank of elements of rings developed in \cite{Sto1}. For a ring (or algebra) $R$ and an element $a \in R$, we study the structure of the ring $aRa$, where $a$ satisfies some finite rank condition.

If $e$ is an idempotent, then the ring $eRe$ is usually called the \emph{(Peirce) corner ring} of $R$ with respect to $e$, so we shall call $aRa$ the \emph{generalized (Peirce) corner ring} of $R$ with respect to $a$.
Corner rings frequently come into play in the structure theory of associative rings, where they often take the role of the building blocks for bigger rings (see \cite{Lam}).
They are also extremely important in Morita theory of equivalences (see \cite{Lam2}, in particular, Corollaries~18.35 and 18.37) and often appear in considerations in several other areas of ring theory, such as extensions of rings, Boolean algebras, rings of operators, path algebras of quivers, etc.
Even in the context of more general type of corner rings introduced in \cite{Lam3}, Peirce corner rings are an example of a kind of a `prototype' for corner rings.
With such a wide range of applications, Peirce corner rings certainly deserve additional attention.

Throughout this paper we will be using the notion of \emph{rank} of an element of a ring, so we recall the definition. For the background and properties of rank we refer the reader to \cite{Sto1}.

\begin{definition}
An element $a \in R$ has \emph{right rank} $0$ if and only if $a=0$. An element $a \in R$ has \emph{right rank} $1$ if and only if $a \neq 0$ and $a$ is contained in some minimal right ideal of $R$. An element $a \in R$ has \emph{right rank} $n>1$ if and only if $a$ is contained in a sum of $n$ minimal right ideals of $R$, but is not contained in any sum of less than $n$ minimal right ideals of $R$. An element $a \in R$ has \emph{infinite right rank} if and only if $a$ is not contained in any sum of minimal right ideals of $R$. The right rank of $a \in R$ will be denoted by $\rank_r a$. The left rank of an element $a \in R$ is defined analogously and denoted by $\rank_l a$.
\end{definition}

In this paper we first describe the structure of corner rings $eRe$, where $e$ is an idempotent of finite rank in a ring $R$ (see Theorem~\ref{semisimple_corner}). The rank $1$ case, is a well known result from the theory of idempotents, which states that a rank one idempotent $e$ gives rise to a division ring $eRe$ (see Proposition~\ref{rank_one}). Our theorem is thus a generalization of this result to arbitrary finite rank. It turns out that the rank of $e$ is a kind of a measure for the size of $eRe$.

In the second part we describe the structure of generalized corner rings $aRa$, where $a$ need not be an idempotent. A partial description of the structure in the setting of semisimple Banach algebras has been given by Bre\v sar and \v Semrl in \cite[Main Theorem~(F)]{Bre-Sem},
however their description does not produce a direct sum decomposition of $aRa$, but only an orthogonal sum decomposition (a decomposition as a sum of additive groups, such that any two of them multiply to $0$). Our main result, Theorem~\ref{main_corner} (along with Corollary~\ref{finite_corner}), gives a direct sum decomposition of $aRa$ for any regular element $a$ of an algebra $R$, such that $a^2$ is regular and has finite rank. At the end we give an example to show that the condition on the rank of $a^2$ is necessary and optimal.

\section{Preliminaries}

All rings and algebras considered in this paper will be associative and unital. For a unital ring $R$ we denote by $M_n(R)$ the ring of all $n \times n$ matrices with entries in $R$. More generally, the set of all $n \times m$ matrices with entries in $R$ will be denoted by $M_{n,m}(R)$.
Standard matrix units in $M_n(R)$ will be denoted by $E_{ij}$, $1 \leq i,j \leq n$, thus $E_{ij}$ is a matrix whose only nonzero entry is entry $(i,j)$ and is equal to $1$. We will also need the ring of all upper triangular $n \times n$ matrices over $R$, denoted by $T_n(R)$.

For the group of all multiplicatively invertible elements of a ring $R$ we will use the standard notation $\mathcal{U}(R)$.
Recall that an element $a \in R$ is called \emph{regular} if there exists an element $b \in R$ such that $a=aba$. In this case $ab$ and $ba$ are idempotents.
If there exists $b \in \mathcal{U}(R)$ that satisfies this condition, then the element $a$ is called \emph{unit-regular}. Equivalently, $a$ is unit-regular if there exists $x \in \mathcal{U}(R)$ and an idempotent $e \in R$, such that $a=ex$.

The \emph{right socle} of a ring $R$ is defined as the sum of all minimal right ideals of $R$. In other words, this is just the set of all element of finite right rank. The \emph{left socle} of $R$ is defined analogously via left ideals. In a semiprime ring $R$, the left and the right socle coincide and are thus simply called the \emph{socle} of $R$ and denoted by $\soc R$.

\section{The idempotent case}

In this section we discuss the structure of the corner ring $eRe$, where $e \in R$ is an idempotent of finite (right) rank. If $e$ has right rank $1$ (i.e. $eR$ is a minimal right ideal), the result is well known, namely

\begin{proposition}[{\cite[Proposition~21.16]{Lam}}]\label{rank_one}
Let $R$ be a ring and $e \in R$ an idempotent. If $eR$ is a minimal right ideal of $R$, then $eRe$ is a division ring. The converse holds if $R$ is a semiprime ring.
\end{proposition}

To generalize this to arbitrary finite rank, we need a technical lemma.

\begin{lemma}\label{corner}
Let $e \in R$ be an idempotent. If $K$ is a minimal right ideal of $R$, then $eKe$ is either zero or a minimal right ideal of $eRe$.
\end{lemma}

\proof By the right-hand side version of \cite[Theorem~21.11]{Lam}, $J \mapsto JR$ defines an injective inclusion-preserving map from the set of right ideals of $eRe$ to the set of right ideals of $R$. Since the right ideal $eKe$ of $eRe$ maps to $eKeR \subseteq K$ and $K$ is a minimal right ideal of $R$, $eKe$ is either zero or a minimal right ideal of $eRe$. \endproof

We say that $a=a_1+a_2+\ldots+a_n$ is a \emph{minimal right decomposition} of $a$, if all $a_i$ have right rank $1$ and $n$ is the right rank of $a$ (see \cite{Sto1}). Next theorem generalizes Proposition~\ref{rank_one}.

\begin{theorem}\label{semisimple_corner}
Let $e \in R$ be an idempotent. If $e$ has finite right rank $n$, then $eRe \cong M_{n_1}(D_1) \times M_{n_2}(D_2) \times \ldots \times M_{n_k}(D_k)$ for some division rings $D_i$ and positive integers $n_i$, where $n_1+n_2+\ldots+n_k=n$. The converse holds if $R$ is a semiprime ring.
\end{theorem}

\proof Let $e$ be an idempotent with finite right rank $n$ and let $e=e_1+e_2+\ldots+e_n$ be its minimal right decomposition. By \cite[Proposition~4.3]{Sto1},
the elements $e_i$ are pairwise orthogonal idempotents of right rank $1$. In particular, $e_iR$ is a minimal right ideal of $R$. Lemma~\ref{corner} implies that $e_iRe=e(e_iR)e$ is a minimal right ideal of $eRe$, hence $e_i$ has right rank $1$ even when considered as an element of $eRe$.
Thus, by \cite[Corollary~4.4]{Sto1},
the right rank of $e$, when considered as an element of $eRe$, is $n$ as well. Since $eRe=e_1Re+e_2Re+\ldots+e_nRe$, and each $e_iRe$ is a simple right $eRe$-module, $eRe$ is a semisimple ring.
Thus, by Wedderburn-Artin Theorem, it is of the desired form. To prove that $n_1+n_2+\ldots+n_k=n$, observe that $e \in eRe$ is a sum of $n_1+n_2+\ldots+n_k$ pairwise orthogonal idempotents, each of which is some standard matrix unit $E_{ii}$ in some $M_{n_j}(D_j)$.
Hence, by \cite[Corollary~4.4]{Sto1},
its right rank in $eRe$ is $n_1+n_2+\ldots+n_k$. Combining this with the above observations, we conclude that $n_1+n_2+\ldots+n_k=n$.

Now let $R$ be a semiprime ring and $eRe$ as in the theorem. As above, $e \in eRe$ is a sum of $n_1+n_2+\ldots+n_k$ orthogonal idempotents, each of which is some standard matrix unit $E_{ii}$ in some $M_{n_j}(D_j)$. So if $ere$ is one of these idempotents, then $(ere)R(ere)=(ere)(eRe)(ere) \cong D_j$ for some $j$. Since $R$ is semiprime, Proposition~\ref{rank_one} implies that $ere$ has right rank $1$ in $R$. Therefore, by \cite[Corollary~4.4]{Sto1},
the idempotent $e$ has right rank $n_1+n_2+\ldots+n_k$ in $R$. \endproof

Theorem~\ref{semisimple_corner} in particular implies that $1 \in R$ has finite right rank if and only if $R$ is a finite direct product of full matrix rings over some division rings.

As a corollary to Theorem~\ref{semisimple_corner} we obtain a generalization of \cite[Corollary~3.8]{Sto1}.

\begin{corollary}\label{reg_rank}
If $a \in R$ is a regular element of finite right and finite left rank, then $\rank_r a=\rank_l a$.
\end{corollary}

\proof Let $a$ be an element of finite right and finite left rank. By \cite[Proposition~4.9]{Sto1},
element $a$ is unit-regular, so $a=eu$ for some idempotent $e$ and some $u \in \mathcal{U}(R)$, where the idempotent $e$ has the same left and right rank as $a$ by \cite[Corollary~3.6]{Sto1}.
Theorem~\ref{semisimple_corner} and its left-hand sided version imply that
\[eRe \cong M_{n_1}(D_1) \times \ldots \times M_{n_k}(D_k) \cong M_{m_1}(E_1) \times \ldots \times M_{m_j}(E_j),\]
where $\rank_r e=n_1+\ldots+n_k$ and $\rank_l e=m_1+\ldots+m_j$. By the uniqueness in the Wedderburn-Artin Theorem it follows that $\rank_r e=\rank_l e$. \endproof

We remark that the assumption that both ranks of $a$ are finite is essential (see \cite[Example~3.4]{Sto1}).

Having a matrix representation for the ring $eRe$, a natural question arises, whether the right rank of an element $a$, which is contained in $eRe$, coincides with the right rank of the corresponding $k$-tuple of matrices. In other words, do the ranks of $a \in eRe$ calculated within $R$ and $eRe$ coincide?
This later question is viable even if $e$ has infinite right rank, however, in this case the answer may be negative. For example, in the ring $R=T_2(\mathbb{C})$, the idempotent $E_{11}$ has infinite right rank in $R$, but right rank $1$ in $E_{11}RE_{11} \cong \mathbb{C}$. Nevertheless, the following corollary answers the question in the affirmative if $e$ has finite rank.

\begin{corollary}\label{subring}
Let $e \in R$ be an idempotent and $a \in eRe$. Suppose any of the following conditions holds:
\begin{enumerate}\sep
\item $a$ is regular and has finite right rank in R,
\item $e$ has finite right rank in R.
\end{enumerate}
Then the right ranks of $a$ in $eRe$ and $R$ coincide.
\end{corollary}

\proof Observe that the first part of the corollary implies the second one. Indeed, if $e$ has finite right rank, then $a \in eRe$ has finite right rank as well. In addition, $a$ is regular since, by Theorem~\ref{semisimple_corner}, the ring $eRe$ is isomorphic to a finite direct product of full matrix rings, which is a regular ring.

To prove the first part of the corollary, we first show that the conclusion is true for an idempotent $f \in eRe$. For the purpose of this proof, let $\rank_R a$ and $\rank_{eRe} a$ denote the right rank of $a$ in $R$ and $eRe$ respectively. For $f=0$ there is nothing to prove, so suppose $f \neq 0$. Choose a minimal right decomposition of $f$ in $R$, say $f=f_1+f_2+\ldots+f_n$, where $n=\rank_R f$.
Observe that $f=efe=ef_1e+ef_2e+\ldots+ef_ne$ is another minimal right decomposition of $f$ in $R$. Hence, $ef_ieR$ is a minimal right ideal of $R$, and by Lemma~\ref{corner}, $ef_ieRe$ is a minimal right ideal of $eRe$. This shows that $\rank_{eRe} ef_ie=1$. From \cite[Proposition~4.3 and Corollary~4.4]{Sto1}
we conclude that $\rank_{eRe} f=n=\rank_R f$.
Now let $a \in eRe$ be regular in $R$. Then $a$ is regular as an element of $eRe$ as well. The same argument as for $f$ shows that $\rank_{eRe} a \leq \rank_R a$. So as an element of $eRe$, $a$ is regular and has finite right rank. By \cite[Proposition~4.9]{Sto1},
$a$ is unit-regular in $eRe$, so there exists an idempotent $h \in eRe$ and $x \in \mathcal{U}(eRe)$ such that $a=hx$.
Let $y$ denote the inverse of $x$ in $eRe$, so that $xy=e$. We have $\rank_{eRe} a=\rank_{eRe} h$ and due to $h=he=hxy=ay$ we also have $\rank_R h=\rank_R ay \leq \rank_R a =\rank_R hx \leq \rank_R h$, i.e. $\rank_R a=\rank_R h$. We have already seen that the two ranks of an idempotent $h$ are the same, thus $\rank_{eRe} a=\rank_R a$. \endproof

\section{The general case}

Now we consider generalized corner ring $aRa$, where $a$ need not be an idempotent. If $e$ is an idempotent of finite (right) rank, then by Theorem~\ref{semisimple_corner}, the Jacobson radical of $eRe$ is zero. For a non-idempotent element $a$ the Jacobson radical of $aRa$ may be very big. In fact, if $a^2=0$, then $aRa$ is a ring with trivial multiplication and its Jacobson radical is the whole ring $aRa$.
So first we want to describe the Jacobson radical of the ring $aRa$.
We shall denote the Jacobson radical of a ring $R$ by $J(R)$. Recall that for an arbitrary idempotent $e \in R$ the Jacobson radical of $eRe$ is equal to $J(eRe)=eRe \cap J(R)=eJ(R)e$ (see \cite{Lam}).

Let $(R,+,\cdot)$ be a ring and $s$ an element of $R$. Define a new multiplication on $R$ by
\[x\ast_s y=xsy\]
for all $x,y \in R$. Then $(R,+,\ast_s)$ is again a ring, which we will denote by $R_s$. We remark that the ring $R_s$ is unital if and only if $s \in \mathcal{U}(R)$. If $R$ is an $F$-algebra, then $R_s$ is also an $F$-algebra for the same multiplication by scalars.

\begin{proposition}\label{deform_1}
Let $a \in R$ be a regular element, $b \in R$ an element, such that $a=aba$, and let $e=ab$. Then $aRa$ is isomorphic to $(eRe)_{eae}$. The isomorphism is given by $x \mapsto xb$ and its inverse is given by $x \mapsto xa$.
\end{proposition}

\proof Let $f: aRa \to eRe$ and $g: eRe \to aRa$ be maps defined by $f(x)=xb$ and $g(x)=xa$. For every $r \in R$ we have $(ara)b=abarab=eare \in eRe$ and $(ere)a=abraba=abra \in aRa$, so the maps $f$ and $g$ are well defined. Observe that $aba=a$ implies $xba=x=abx$ for all $x \in aRa$. In addition, $xab=xe=x$ for all $x \in eRe$. This shows that $g$ is the inverse of $f$. Clearly the map $f$ is additive (linear, if $R$ is an algebra), so it remains to show that it is also multiplicative as a map from $aRa$ to $(eRe)_{eae}$. Let $x,y \in aRa$ be arbitrary. By the above we have
\[f(x)\ast_{eae} f(y)=(xb)(eae)(yb)=(xb)(a^2b)(yb)=(xba)(aby)b=xyb=f(xy),\]
as required. \endproof

\begin{proposition}\label{deform_2}
For every $s \in R$ we have
\[J(R_s)=\set{x \in R}{sxs \in J(R)}.\]
\end{proposition}

\proof Recall that the Jacobson radical of a ring can be characterized as the set of all elements $x$, such that $xr$ is right quasi-regular for every element $r$.
An element $a$ is right quasi-regular if there exists some element $b$ such that $a+b-ab=0$.

Suppose $x \in J(R_s)$ and choose an arbitrary $r \in R$. Then $x \ast_s r$ is right quasi-regular in $R_s$, so there exists $y \in R$ such that $x\ast_s r+y-x\ast_s r \ast_s y=0$, or equivalently $xsr+y-xsrsy=0$. Multiplying from the left by $s$ we get $(sxs)r+(sy)-(sxs)r(sy)=0$. This shows that $(sxs)r$ is right quasi-regular in $R$ for every $r \in R$, thus $sxs \in J(R)$.

Now suppose $sxs \in J(R)$ and choose an arbitrary $r \in R$. Then there exists $y \in R$ such that $(sxs)r+y-(sxs)ry=0$. This implies $y=sz$ where $z=xsry-xsr$, so that $(sxs)r+sz-(sxs)rsz=0$ or equivalently $s(xsr+z-xsrsz)=0$. Now denote $w=xsr+z-xsrsz$, so that $sw=0$. These two equalities imply $xsr+(z-w)-xsrs(z-w)=0$ or equivalently $x\ast_s r+(z-w)-x\ast_s r\ast_s (z-w)=0$. Hence $x\ast_s r$ is right quasi-regular in $R_s$. Since $r$ was arbitrary, we conclude that $x \in J(R_s)$. \endproof

Proposition~\ref{deform_2} in particular implies $J(R) \subseteq J(R_s)$. We can now describe the Jacobson radical of $aRa$.

\begin{corollary}\label{jara}
For a regular element $a \in R$ we have
\[J(aRa)=\set{x \in aRa}{axa \in J(R)}.\]
\end{corollary}

\proof Choose $b \in R$ such that $aba=a$ and denote $e=ab$. Propositions~\ref{deform_1} and \ref{deform_2} imply
\[J(aRa)=\set{x \in aRa}{xb \in J((eRe)_{eae})}=\set{x \in aRa}{a^2bxba^2b \in J(eRe)}.\]
For every $x \in aRa$ we have $a^2bxba^2b=a(abxba)ab=axab=axe$, hence
\[J(aRa)=\set{x \in aRa}{axe \in J(eRe)}=\set{x \in aRa}{axe \in eJ(R)e}.\]
It suffices to prove that $axe \in eJ(R)e$ is equivalent to $axa \in J(R)$. Suppose $axe \in eJ(R)e$. Multiplying from the right by $a$ we obtain $axa \in eJ(R)a \subseteq J(R)$. Now suppose $axa \in J(R)$. Multiplying from the left by $e$ and from the right by $be$ we get $axe \in eJ(R)be \subseteq eJ(R)e$. \endproof

Corollary~\ref{jara} states that for a regular element $a \in R$, $J(aRa)$ is the largest subset of $R$ that satisfies
\[aJ(aRa)a=a^2Ra^2 \cap J(R).\]

\begin{lemma}\label{lem1}
Let $s,u,v$, and $e$ be elements of $R$, where $u$ and $v$ are invertible and $e$ is an idempotent. Then
\begin{enumerate}
\item $R_{usv} \cong R_s$, where the isomorphism is given by $x \mapsto vxu$,
\item\label{c2} $\displaystyle R_e \cong \left[\begin{array}{ccc} 0 & (1-e)Re & (1-e)R(1-e) \\  0 & eRe & eR(1-e) \\ 0 & 0 & 0 \\ \end{array}\right]$, where the isomorphism is induced by the Peirce decomposition,
\item\label{c3} $R_e/J(R_e) \cong eRe/J(eRe)$, where the isomorphism is induced by the map $x \mapsto exe$.
\end{enumerate}
\end{lemma}

\proof Let $f: R_{usv} \to R_s$ be a map defined by $f(x)=vxu$. Clearly, $f$ is an additive bijection (linear, if $R$ is an algebra). Since $f(x \ast_{usv} y)=f(xusvy)=vxusvyu=f(x)sf(y)=f(x) \ast_s f(y)$ for all $x,y \in R$, the map $f$ is a ring isomorphism.

The isomorphism in \ref{c2} is given by
\[g: x \mapsto \left[\begin{array}{ccc} 0 & (1-e)xe & (1-e)x(1-e) \\  0 & exe & ex(1-e) \\ 0 & 0 & 0 \\ \end{array}\right].\]
By means of Peirce decomposition its inverse is easily seen to be just the summation of elements of the matrix. Clearly, $g$ is additive (linear, if $R$ is an algebra). An easy computation shows that $g$ is also multiplicative.

\ref{c3} is an easy consequence of \ref{c2}. \endproof

We can now describe the structure of the ring $aRa$ in terms of its Jacobson radical. Recall that by \cite[Corollary~3.8]{Sto1}
the left and right rank of any element in a semiprime ring coincide, so we omit the adjectives in this case and simply speak of rank.

\begin{theorem}\label{a_corner}
Let $R$ be a unital semiprime ring and $a \in R$ an element of finite rank. Then $J(aRa)=\set{x \in aRa}{axa=0}$ and $(aRa)/J(aRa) \cong M_{n_1}(D_1) \times M_{n_2}(D_2) \times \ldots \times M_{n_k}(D_k)$ for some division rings $D_i$ and positive integers $n_i$, where $n_1+n_2+\ldots+n_k=\rank a^2$.
\end{theorem}

\proof By \cite[Theorem~4.10]{Sto1},
element $a$ is unit-regular, hence $a=eu$ for some $u \in \mathcal{U}(R)$ and some idempotent $e \in R$ of finite rank.
Let $b=u^{-1}$, so that $aba=a$ and $e=ab$. Corollary~\ref{jara} implies that $J(aRa)=\set{x \in aRa}{axa \in J(R)}$.
Observe that $(\soc R \cap J(R))^2=0$, since $J(R)$ annihilates all simple $R$-modules. The semiprimeness of $R$ therefore implies $(\soc R \cap J(R))=0$. But $axa \in \soc R$ for any $x \in R$, so $J(aRa)=\set{x \in aRa}{axa=0}$.
Applying Proposition~\ref{deform_1} we see that $aRa \cong (eRe)_{eae}$. Theorem~\ref{semisimple_corner} in particular implies that $eae$ is unit-regular in $eRe$. Let $eae=fw$, where $f$ is an idempotent in $eRe$ and $w \in \mathcal{U}(eRe)$.
Lemma~\ref{lem1} implies $(eRe)_{fw} \cong (eRe)_f$ and $(eRe)_f/J((eRe)_f) \cong f(eRe)f/J(f(eRe)f) = fRf/J(fRf)$. Since $f \in eRe$ has finite rank as well, Theorem~\ref{semisimple_corner} implies that the ring $fRf$ has zero Jacobson radical and is isomorphic to $M_{n_1}(D_1) \times M_{n_2}(D_2) \times \ldots \times M_{n_k}(D_k)$ for some division rings $D_i$, where $n_1+n_2+\ldots+n_k=\rank f$.
Putting everything together we see that $aRa/J(aRa) \cong fRf$, so the conclusion of the theorem will follow as soon as we prove that $\rank f = \rank a^2$. Invoking to Corollary~\ref{subring} we infer that the rank of $f$ is equal to the rank of $eae=a^2u^{-1}$, which is further equal to the rank of $a^2$. \endproof

Let $R$ and $a$ be as in Theorem~\ref{a_corner}. The proof of the theorem shows that there exist idempotents $e$ and $f$, and invertible elements $u$ and $v$, such that $a=eu$, $a^2=fv$ and $f=efe$ (the existence of $v$ follows from \cite[Corollary~4.7]{Sto1}).
In this case the isomorphism $aRa/J(aRa) \cong fRf$ is induced by the map $x \mapsto axu^{-1}f=axav^{-1}$, as can be seen by examining the proof carefully.

As a corollary we can characterize when $aRa$ is a semisimple ring.

\begin{corollary}
Let $R$ be a unital semiprime ring and $a$ an element of finite rank. The following conditions are equivalent:
\begin{enumerate}\sep
\item $aRa$ is a semiprime ring,
\item $\rank a^2=\rank a$,
\item there exist $b,c \in R$ such that $a=a^2b=ca^2$,
\item $J(aRa)=0$,
\item $aRa \cong M_{n_1}(D_1) \times M_{n_2}(D_2) \times \ldots \times M_{n_k}(D_k)$ for some division rings $D_i$ and positive integers $n_i$, where $n_1+n_2+\ldots+n_k=\rank a$.
\end{enumerate}
\end{corollary}

\proof Suppose $aRa$ is a semiprime ring, but $\rank a^2 <\rank a$. By \cite[Theorem~4.10]{Sto1}
there exist finite rank idempotents $e$ and $f$, with $\rank f <\rank e$, and invertible elements $u$ and $v$, such that $a=eu$ and $a^2=fv$. This in particular implies $efv=ea^2=a^2=fv$, hence $ef=f$. Observe that $e \neq fe$, since the ranks of the two are different.
Thus $(1-f)e \neq 0$ and the semiprimeness of $R$ implies the existence of an element $t$ such that $(1-f)et(1-f)e \neq 0$. Let $r=u^{-1}t(1-f)$. Then $(1-f)arau^{-1}=(1-f)(eu)u^{-1}t(1-f)(eu)u^{-1}=(1-f)et(1-f)e \neq 0$, so that $ara \neq 0$. However, $ra^2=u^{-1}t(1-f)fv=0$, so that $(ara)(aRa)(ara)=0$. This contradicts the assumption that $aRa$ is a semiprime ring.

Assume $\rank a^2 =\rank a$. Then by \cite[Corollary~4.11]{Sto1}
there exists an invertible element $x$ such that $a^2=ax$. Hence $a=a^2x^{-1}$, so we may take $b=x^{-1}$. The existence of $c$ is proved similarly, because the rank in semiprime rings is left-right symmetric.

Now let $b$ and $c$ be elements such that $a=a^2b=ca^2$. By \cite[Theorem~4.10]{Sto1},
$a$ is a regular element, hence $J(aRa)=\set{x \in aRa}{axa \in J(R)}$ by Corollary~\ref{jara}. So if $x \in J(aRa)$, then $axa \in J(R)$. Multiplying from the left by $c$ and from the right by $b$, and taking into account that $x \in aRa$, we get $x \in cJ(R)b \subseteq J(R)$. Consequently, $x \in \soc R \cap J(R)=0$, since $R$ is semiprime. Hence $J(aRa)=0$.

Theorem~\ref{a_corner} shows that $(iv)$ implies $aRa \cong M_{n_1}(D_1) \times M_{n_2}(D_2) \times \ldots \times M_{n_k}(D_k)$, where $n_1+n_2+\ldots+n_k=\rank a^2$. But this clearly implies $(i)$, which implies $(ii)$, so $n_1+n_2+\ldots+n_k=\rank a$. \endproof

For an element $r \in R$, let $(r)$ denote the ideal of $R$ generated by $r$. The following lemma will be used in the proof of our main theorem.

\begin{lemma}\label{cap_dot}
Let $f \in R$ be an idempotent of finite right rank and $K$ a right ideal of $R$. Then $(f)\cap K=K \cdot (f)$.
\end{lemma}

\proof Clearly $K \cdot (f) \subseteq (f)\cap K$, so suppose $(f)\cap K \not\subseteq K \cdot (f)$. Let $f=f_1+f_2+\ldots+f_n$ be a minimal right decomposition of $f$. By \cite[Proposition~4.3]{Sto1},
$f_i$ are orthogonal idempotents of right rank $1$. In particular, $f_i \in (f)$. Every element in $(f)$ is of the form $\sum_{i=1}^k x_ie_iy_i$ for some $e_i \in \{f_1,f_2,\ldots,f_n\}$ and $x_i,y_i \in R$. Let $k$ be the least positive integer, such that $s=\sum_{i=1}^k x_ie_iy_i \in (f)\cap K \backslash K \cdot (f)$ for some $e_i \in \{f_1,f_2,\ldots,f_n\}$ and $x_i,y_i \in R$. Since $e_ky_k \neq 0$ and $e_k$ has right rank $1$, there exists $z \in R$ such that $e_ky_kz=e_k$. Let $g=ze_ky_k \in (f)$, so that $e_ky_kg=e_ky_k$. Then
$sg=\sum_{i=1}^{k-1} x_ie_iy_ig +x_ke_ky_k$ and hence $s-sg=\sum_{i=1}^{k-1} x_ie_iy_i(1-g)$. Clearly $s-sg=s(1-g) \in (f)\cap K$, so the choice of $k$ implies $s-sg \in K \cdot (f)$ (this is true even if $k=1$, in which case $s-sg=0$). Since $sg \in K \cdot (f)$, we get a contradiction $s \in K \cdot (f)$. \endproof

Building off of Theorem~\ref{a_corner} we now describe the structure of $aRa$ in terms of its ideals. As indicated by Theorem~\ref{a_corner}, it is crucial that $a^2$ has finite rank, if we want to describe the structure of $aRa$. The rank of $a$ does not seem to play much of a role in the matter.
So we will only assume that the rank of $a^2$ is finite, while the rank of $a$ may be infinite. We will also not assume $R$ to be semiprime, instead, we will work with regular elements. In addition, we will need $R$ to be an algebra over a field. Compare next theorem and its corollary with \cite[Lemma~2.7 and Main Theorem~(F)]{Bre-Sem}.

\begin{theorem}\label{main_corner}
Let $F$ be a field, $R$ a unital $F$-algebra, and $a \in R$ a regular element, such that $a^2$ is regular as well. If $a^2$ has finite right rank $n$, then there exist ideals $I_0,I_1,\ldots,I_k \ideal aRa$, such that
\begin{enumerate}
\item\label{mc1} $aRa = I_0 \oplus \bigoplus_{j=1}^k I_j$,
\item\label{mc2} $I_0^2=0$, while $I_1,I_2,\ldots,I_k$ are directly irreducible and $I_j/J(I_j) \cong M_{n_j}(D_j)$, for some division algebras $D_j$ and positive integers $n_j$,
\item\label{mc3} $J(I_j)^3=0$ for all $1 \leq j \leq k$,
\item\label{mc4} $n_1+n_2+\ldots+n_k=n$.
\end{enumerate}
The converse holds if $R$ is a semiprime algebra.

Up to permutation the ideals $I_1,I_2,\ldots,I_k$ are uniquely determined by \ref{mc1} and \ref{mc2}, while $I_0$ is unique only up to isomorphism. Moreover, $I_1,I_2,\ldots,I_k$ are principal ideals generated by any nonzero idempotent they contain.

If $R$ is a prime algebra, then $k\leq 1$.
\end{theorem}

\proof The beginning of the proof is similar to that of Theorem~\ref{a_corner}, just a lot more care is needed, because the right rank of $a$ need not be finite. Choose $b,c \in R$ such that $aba=a$ and $a^2ca^2=a^2$, and let $e=ab$, which is an idempotent, possibly of infinite right rank.
By Proposition~\ref{deform_1}, $aRa \cong (eRe)_{eae}$.
Observe that $eae=a^2b$ is a regular element of $eRe$. This is because $(a^2b)(acab)(a^2b)=a^2b$ and $acab=eace \in eRe$. In addition, $eae=a^2b$ has finite right rank in $R$, since $a^2$ has finite right rank.
By Corollary~\ref{subring}, the right ranks of $eae$ in $R$ and $eRe$ coincide. Hence, by \cite[Proposition~4.9]{Sto1},
$eae$ is unit-regular in $eRe$.
Let $eae=fw$, where $f=efe$ is an idempotent in $eRe$ and $w$ is an invertible element of $eRe$. Then $f=eaew^{-1}$, where the inverse is taken in $eRe$. This shows that the right rank of $f$ in $R$ is finite. By Corollary~\ref{subring}, the right ranks of $f$ in $R$ and $eRe$ coincide.
In $eRe$ the right rank of $f$ is the same as the right rank of $eae$, hence the same holds in $R$. But, by \cite[Proposition~3.5]{Sto1},
the right rank of $eae=a^2b$ in $R$ is the same as the right rank of $a^2$, since $a^2=eaea$.
We conclude that $\rank_r f=\rank_r a^2$. As in the proof of Theorem~\ref{a_corner}, $aRa \cong (eRe)_f$ and $fRf \cong M_{n_1}(D_1) \times M_{n_2}(D_2) \times \ldots \times M_{n_k}(D_k)$ for some division algebras $D_j$, where $n_1+n_2+\ldots+n_k=\rank_r a^2=n$.

We may work in $(eRe)_f$ since everything carries over to $aRa$. Let $f_j$ be the idempotent in $fRf$ that corresponds under the above isomorphism to the identity matrix in $M_{n_j}(D_j)$, so that $f_1,f_2,\ldots,f_k$ are orthogonal idempotents (in $eRe$ as well as in $(eRe)_f$) with $f=f_1+f_2+\ldots+f_k$. This in particular implies that $f_j=ff_j$ has finite right rank in $R$ (but not necessarily in $(eRe)_f$). Note that $f_jRf_j=f_j(fRf)f_j \cong M_{n_j}(D_j)$ for all $1\leq j\leq k$ and $f_iRf_j=f_i(fRf)f_j=0$ for all $i \neq j$.
Define $f_0=e-f$ and observe that $f_0$ is an idempotent in $R$ orthogonal to $f$ and hence to all $f_j$, $1 \leq j \leq k$. Now let $I_j$, $1 \leq j \leq k$, be the ideal of $(eRe)_f$ generated by $f_j$. Since $e=f_0+f_1+\ldots+f_k$ and $f_iRf_j=0$ for $i \neq j$, $1\leq i,j\leq k$, we have
\begin{align}
I_j &=(eRe)\ast_f f_j \ast_f (eRe)=eRf_jRe= \nonumber\\
&= f_jRf_jRf_j+f_jRf_jRf_0+f_0Rf_jRf_j+f_0Rf_jRf_0= \label{I_j}\\
&= f_jRf_j+f_jRf_0+f_0Rf_j+f_0Rf_jRf_0 \nonumber.
\end{align}
Denote $N_j=f_jRf_0+f_0Rf_j+f_0Rf_jRf_0$. The orthogonality of $f$ and $f_0$ implies $N_j \ast_f I_j \subseteq f_0Rf_j+f_0Rf_jRf_0$, $I_j \ast_f N_j \subseteq f_jRf_0+f_0Rf_jRf_0$, and $I_j \ast_f N_j \ast_f I_j=0$, so $N_j$ is a nilpotent ideal of $I_j$ of nilindex $\leq 3$. On the other hand, $f_jRf_j\cong M_{n_j}(D_j)$ as a subalgebra of $R$ and as a subalgebra of $(eRe)_f$. This shows that $J(I_j)=N_j$ and $I_j/J(I_j) \cong M_{n_j}(D_j)$. Suppose $I_j=K \oplus L$, where $K$ and $L$ are ideals of $I_j$. Then $I_j/J(I_j) \cong K/J(K) \oplus L/J(L)$. This implies $K/J(K)=0$ or $L/J(L)=0$, since $I_j/J(I_j)$ is a simple algebra. We may assume $K/J(K)=0$, thus $K \subseteq J(K) \subseteq J(I_j)=N_j$. Now write $f_j=p+l$, where $p \in K \subseteq N_j$ and $l \in L$. Then by the above $f_j=f_j \ast_f p \ast_f f_j+f_j \ast_f l \ast_f f_j=f_j \ast_f l \ast_f f_j \in L$. Since $I_j$ is generated by $f_j$ it follows that $I_j=L$ and $K=0$. This shows that $I_j$ is directly irreducible.

Since $f_j$, $1 \leq j \leq k$, are orthogonal idempotents and $f=f_1+f_2+\ldots+f_k$, the ideal $\sum_{j=1}^k I_j$ is generated by $f$. Similarly as in \eqref{I_j} we thus have
\begin{equation}\label{direct1}
\sum_{j=1}^k I_j=fRf+fRf_0+f_0Rf+f_0RfRf_0.
\end{equation}

Let $I_0$ be a vector subspace of $f_0Rf_0$, such that $f_0Rf_0=f_0RfRf_0 \oplus I_0$ as vector spaces. Now $I_0 \ast_f (eRe)=I_0fRe\subseteq f_0R(f_0f)Re=0$ and $(eRe) \ast_f I_0=0$, so $I_0$ is in fact a square-zero ideal of $(eRe)_f$.
By the definition of $I_0$,
\[I_0+\sum_{j=1}^k I_j=fRf+fRf_0+f_0Rf+f_0Rf_0=eRe,\]
where the last equality is just the Peirce decomposition of $eRe$. It remains to prove that the sum $I_0+\sum_{j=1}^k I_j$ is direct. The sum of vector spaces on the right-hand side of \eqref{direct1} is direct since the idempotents $f$ and $f_0$ are orthogonal. Together with $I_0 \subseteq f_0Rf_0$ this implies
\[I_0 \cap \sum_{j=1}^k I_j=I_0 \cap f_0RfRf_0=0,\]
where the last equality follows from the definition of $I_0$. For $1 \leq i \leq k$ we have as in \eqref{I_j}
\begin{align*}
I_i \cap \Big(I_0 +\sum_{\substack{j=1 \\ j \neq i}}^k I_j\Big) &=(f_iRf_i+f_iRf_0+f_0Rf_i+f_0Rf_iRf_0) \cap \\
& \cap(f_i'Rf_i'+f_i'Rf_0+f_0Rf_i'+f_0Rf_i'Rf_0+I_0),
\end{align*}
where $f_i'=f-f_i$. Since the idempotents $f_i,f_i'$ and $f_0$ are orthogonal and $I_0 \subseteq f_0Rf_0$, this boils down to
\begin{align*}
I_i \cap \Big(I_0 +\sum_{\substack{j=1 \\ j \neq i}}^k I_j\Big) &=
f_0Rf_iRf_0 \cap (f_0Rf_i'Rf_0+I_0) =f_0Rf_iRf_0 \cap f_0Rf_i'Rf_0,
\end{align*}
where the last equality is a consequence of the definition of $I_0$. Lemma~\ref{cap_dot}, together with $f_iRf_j=0$ for $i \neq j$, $1\leq i,j\leq k$, therefore implies
\[I_i \cap \Big(I_0 +\sum_{\substack{j=1 \\ j \neq i}}^k I_j\Big)
\subseteq Rf_iR \cap Rf_i'R =Rf_i'Rf_iR=0.\]
This shows that the sum $I_0+\sum_{j=1}^k I_j$ is direct.

By definition $I_j$ is generated by $f_j$. Now let $e_j$ be any nonzero idempotent in $I_j$. Since $J(I_j)$ does not contain nonzero idempotents and $I_j/J(I_j)$ is a simple algebra, we have $I_j=(e_j)+J(I_j)$, where $(e_j)$ denotes the ideal of $(eRe)_f$ generated by $e_j$. In particular, $f_j=d_j+h_j$, where $d_j \in (e_j)$ and $h_j \in J(I_j)$. Hence $f_j=f_j^3=(d_j+h_j)^3 \in (e_j)$, because $h_j^3=0$ (all powers here are taken in $(eRe)_f$). This implies $I_j=(e_j)$.

To prove the uniqueness of ideals $I_0,I_1,\ldots,I_k$ suppose ideals $I_0',I_1',\ldots,I_m'$ also satisfy the conditions \ref{mc1} and \ref{mc2}. Let $1 \leq j \leq k$ and write $f_j=f_0'+f_1'+\ldots+f_m'$, where $f_i' \in I_i'$. Since $f_j$ is an idempotent in $(eRe)_f$, the fact that the sum of ideals $I_i'$ is direct implies that $f_i'$ are orthogonal idempotents.
Multiplying the equation by $f_i'$ we infer that $f_i' \in I_j$, hence $I_j$ is generated by $\{f_0', f_1',\ldots,f_m'\}$. Clearly $f_0'=0$, since $I_0'$ is square-zero ideal. Hence $I_j=(f_1')+(f_2')+\ldots+(f_m')$ and this sum is direct because $(f_i') \subseteq I_i'$. As $I_j$ is directly indecomposable, we must have $I_j=(f_i') \subseteq I_i'$ for some $1 \leq i \leq m$. Suppose there is another $j' \neq j$ such that $I_{j'} \subseteq I_i'$. Then $I_j \oplus I_{j'} \ideal I_i'$, hence
\begin{align*}
M_{n_j}(D_j) \oplus M_{n_{j'}}(D_{j'}) &\cong I_j/J(I_j) \oplus I_{j'}/J(I_{j'}) \cong (I_j \oplus I_{j'})/J(I_j \oplus I_{j'})=\\
&= (I_j \oplus I_{j'})/(I_j \oplus I_{j'})\cap J(I_i') \cong \\
&\cong (I_j \oplus I_{j'}+J(I_i'))/J(I_i') \ideal \\
&\ideal I_i'/J(I_i') \cong M_{n_i}(D_i),
\end{align*}
which is a contradiction. Now suppose there is $1 \leq i' \leq m$ such that $I_j \not\subseteq I_{i'}'$ for all $1 \leq j \leq k$. Then
\[(I_{i'}')^2 \subseteq \Big(I_0 \oplus \bigoplus_{j=1}^k I_j\Big)^2 \subseteq I_0^2 \oplus \bigoplus_{j=1}^k I_j^2 \subseteq \bigoplus_{j=1}^k I_j \subseteq \bigoplus_{\substack{i=1 \\ i \neq i'}}^m I_i'.\]
Hence $(I_{i'}')^2=0$, since the sum $\bigoplus_{i=1}^m I_i'$ is direct. This is a contradiction, because $I_{i'}'/J(I_{i'}') \cong M_{n_{i'}}(D_{i'})$.

From all the above we conclude that the inclusion induces a bijection from $\{I_1,I_2,\ldots,I_k\}$ to $\{I_1',I_2',\ldots,I_m'\}$, in particular $m=k$. We may henceforth assume that $I_j \subseteq I_j'$ for all $1 \leq j \leq k$.

Define $S=I_0 \cap \bigoplus_{j=1}^k I_j'$. Choose an arbitrary element $x' \in \bigoplus_{j=1}^k I_j' \subseteq I_0 \oplus \bigoplus_{j=1}^k I_j$ and write it as $x'=x_0+x$, where $x_0 \in I_0$ and $x \in \bigoplus_{j=1}^k I_j \subseteq \bigoplus_{j=1}^k I_j'$. Then $x_0=x'-x \in S$. Hence
\begin{equation}\label{direct2}
\bigoplus_{j=1}^k I_j'=S \oplus \bigoplus_{j=1}^k I_j,
\end{equation}
where this sum is direct because $S \subseteq I_0$. For $1 \leq i \leq k$ define $S_i=I_i' \cap (S \oplus \bigoplus_{j=1, j \neq i}^k I_j)$. Choose an arbitrary $y_i' \in I_i'$. By \eqref{direct2} we can write it as $y_i'=s+\sum_{j=1}^k y_j$, where $s \in S$ and $y_j \in I_j$ for $1 \leq j \leq k$.  Observe that $y_i'-y_i=s+\sum_{j=1,j \neq i}^k y_j \in S_i$, hence $y_i'=(y_i'-y_i)+y_i \in S_i+I_i$. This shows that $I_i'=S_i \oplus I_i$, where this sum is direct because the one in \eqref{direct2} is. Since $I_j'$ is directly irreducible, we conclude that $I_i'=I_i$ for all $1 \leq i \leq k$.

By the above the ideal $I=\bigoplus_{j=1}^k I_j \ideal aRa$ is uniquely determined. Hence, $I_0 \cong aRa/I$ is unique up to isomorphism.

To prove the converse, suppose $R$ is semiprime and \ref{mc1}--\ref{mc4} hold. It suffices to prove that the right rank of $a^2$ is finite, since the first part and the uniqueness of $I_j$, $1 \leq j \leq k$, will then imply that the right rank of $a^2$ is $n$, because $n_j$ is uniquely determined by $I_j$. Fix some $j$, $1 \leq j \leq k$.
Since $J(I_j)$ is a nilpotent ideal, the idempotents in $I_j/J(I_j)$ can be lifted to $I_j$ (see \cite[Theorem~21.28]{Lam}, where the proof is done for unital rings, however, the same can be proved for non-unital rings by simply adjoining a unit to the ring). Let $g_j$ be an idempotent in $I_j$ lifting the identity element of $I_j/J(I_j)$. By \ref{mc1} we have $g_jRg_j=g_jg_jRg_jg_j \subseteq g_j(aRa)g_j=g_jI_jg_j \subseteq g_jRg_j$, hence
\[(g_jRg_j+J(I_j))/J(I_j)=(g_jI_jg_j+J(I_j))/J(I_j)=I_j/J(I_j) \cong M_{n_j}(D_j),\]
because $g_j+J(I_j)$ is the identity element of $I_j/J(I_j)$. On the other hand,
\[(g_jRg_j+J(I_j))/J(I_j) \cong g_jRg_j/(J(I_j)\cap g_jRg_j).\]
By \ref{mc3}, $J(I_j)\cap g_jRg_j$ is a nilpotent ideal of $g_jRg_j$. However, $g_jRg_j$ is a semiprime algebra (since $R$ is), hence $J(I_j)\cap g_jRg_j=0$. All the above now implies $g_jRg_j \cong M_{n_j}(D_j)$, so $g_j$ has finite right rank in $R$ by Theorem~\ref{semisimple_corner}.

Next we prove that the ideal $I_j$ is generated by $g_j$, as this is not part of the assumptions at this point. In what follows $1$ will denote the identity element of $R$. Previous paragraph additionally implies that $I_j=g_jRg_j+J(I_j)$, hence $I_j(1-g_j)=J(I_j)(1-g_j)$. Observe that $J(I_j)(1-g_j)RJ(I_j) \subseteq aRa$, \mbox{$J(I_j)(1-g_j) \subseteq I_j$}, and $(1-g_j)J(I_j) \subseteq I_j$, therefore \ref{mc1} implies
\[J(I_j)(1-g_j)\Big(J(I_j)(1-g_j)RJ(I_j)\Big)(1-g_j)J(I_j) \subseteq J(I_j)(1-g_j)I_j(1-g_j)J(I_j).\]
By the above we have $(1-g_j)I_j(1-g_j)=(1-g_j)J(I_j)(1-g_j) \subseteq J(I_j)$, hence $J(I_j)(1-g_j)I_j(1-g_j)J(I_j)=0$ by \ref{mc3}. Putting everything together we infer $\Big(RJ(I_j)(1-g_j)J(I_j)(1-g_j)R\Big)^2=0$. Since $R$ is a unital semiprime algebra, this implies $I_j(1-g_j)I_j(1-g_j)=J(I_j)(1-g_j)J(I_j)(1-g_j)=0$. Similarly $(1-g_j)I_j(1-g_j)I_j=0$. Define $Z_j=\{z \in I_j\ ; \ I_jz=zI_j=0\}$, which is 
an ideal of $I_j$. We have just proved that $(1-g_j)I_j(1-g_j) \in Z_j$. For every $r \in R$ we have $r=rg_j+g_jr-g_jrg_j+(1-g_j)r(1-g_j)$, therefore $I_j=I_jg_j+g_jI_j+g_jI_jg_j+(1-g_j)I_j(1-g_j)=(g_j)+Z_j$, where $(g_j)$ denotes the ideal of $aRa$ generated by $g_j$. Hence, there is a vector subspace $V_j \subseteq Z_j$, such that $I_j=(g_j) \oplus V_j$ as vector spaces.
However, by definition of $Z_j$, any subspace of $Z_j$ is clearly an ideal of $I_j$, so the above direct sum is a direct sum of ideals. By \ref{mc2}, $I_j$ is directly indecomposable, therefore $V_j=0$ and $I_j$ is generated by $g_j$. This in particular implies, that every element in $I_j$ has finite right rank in $R$. So it suffices to prove that $a^2 \in \bigoplus_{j=1}^k I_j$.
Recall that due to the regularity of $a^2$ we have $a^2ca^2=a^2$. Hence \ref{mc1} and \ref{mc2} imply $a^2=a^2ca^2ca^2ca^2=(a^2ca)(aca)(aca^2) \in (aRa)^3 \subseteq \bigoplus_{j=1}^k I_j$, as required.

Finally, if $R$ is a prime algebra, then $fRf \cong M_{n_1}(D_1) \times M_{n_2}(D_2) \times \ldots \times M_{n_k}(D_k)$ is a prime algebra as well, so that either $fRf=0$ or $k=1$. \endproof

In the theory of associative rings, in order to prove a structure theorem for certain kind of rings, one usually has to either assume some kind of finiteness condition on one-sided ideals of the ring or work with idempotents with special properties, e.g. central idempotents.
In our situation, in Theorem~\ref{main_corner}, the role of the finiteness condition is taken by the finite rank condition. However, since we only assume $a^2$ to have finite rank, while $a$ may have infinite rank, the ring $aRa$ still retains certain aspects of rings without finiteness conditions. In particular, while the factor ring $aRa/J(aRa)$ is right artinian, the ring $aRa$ need not be.
Observe also, that in general the idempotents $f_j$, that induce the decomposition in Theorem~\ref{main_corner}, are not central.

An important point of Theorem~\ref{main_corner} is that, for an element $a$ with square of finite rank, the direct sum decomposition of the factor ring $aRa/J(aRa)$ induces a direct sum decomposition of the whole ring $aRa$. This is somewhat similar to the situation for commutative artinian rings (cf. \cite[Corollary~2.16]{Eis}). For more general noncommutative rings, e.g. upper triangular matrices $T_n(F)$, this is not true.

\begin{corollary}\label{finite_corner}
Suppose that the element $a$ from Theorem~\ref{main_corner} has finite right rank. Then
\begin{enumerate}
\item 
$\displaystyle I_j \cong \left[\begin{array}{ccc} 0 & M_{m_j,n_j}(D_j) & M_{m_j}(D_j) \\  0 & M_{n_j}(D_j) & M_{n_j,m_j}(D_j) \\ 0 & 0 & 0 \\ \end{array}\right]$\\
for some nonnegative integer $m_j$, for all $j=1,2,\ldots,k$,
\item 
$\displaystyle I_0 \cong \left[\begin{array}{cc} 0 & S \\ 0 & 0 \\ \end{array}\right]$, where
$\displaystyle S=\bigoplus_{j=k+1}^l M_{m_j}(D_j)$\\
for some nonnegative integers $m_{k+1},m_{k+2},\ldots,m_l$ and division algebras $D_{k+1},D_{k+2},\ldots,D_l$,
\item $m_1+m_2+\ldots+m_l+n_1+n_2+\ldots+n_k=\rank_r a$.
\end{enumerate}
\end{corollary}

\proof Assume all notations are as in the proof of Theorem~\ref{main_corner}. Idempotent $e$ has finite rank in this case. Theorem~\ref{semisimple_corner} implies $eRe \cong M_{p_1}(E_1) \times M_{p_2}(E_2) \times \ldots \times M_{p_l}(E_l)$ for some division algebras $E_j$. We may write idempotent $f$ as $f=h_1+h_2+\ldots+h_l$, where $h_j \in M_{p_j}(E_j)$ are idempotents as well. Some of $h_j$ may be zero. By rearranging and renumbering the terms in the decomposition, we may assume that $h_j \neq 0$ for $j\leq k'$ and $h_j=0$ for $j>k'$. Clearly, $aRa \cong (eRe)_f \cong M_{p_1}(E_1)_{h_1} \times M_{p_2}(E_2)_{h_2} \times \ldots \times M_{p_l}(E_l)_{h_l}$. Let $J_j$, $j=0,1,2,\ldots,k'$, be ideals of $aRa$ such that $J_0 \cong \prod_{j=k'+1}^l M_{p_j}(E_j)_{h_j}$ and $J_j \cong M_{p_j}(E_j)_{h_j}$ for $1\leq j\leq k'$.
For $j>k'$ the multiplication in $M_{p_j}(E_j)_{h_j}$ is trivial, so clearly $J_0 \cong \left[\begin{array}{cc} 0 & S \\ 0 & 0 \\ \end{array}\right]$, where $S=\bigoplus_{j=k'+1}^l M_{p_j}(E_j)$. Let $1\leq j\leq k'$. The reduced column echelon form of the reduced row echelon form of $h_j$ is a diagonal idempotent $d_j$ with all the entries which are equal to $1$ collected at the beginning of the diagonal. Hence $h_j=u_jd_jv_j$ for some invertible elements $u_j$ and $v_j$. Let $n_j'$ be the matrix rank of $d_j$ and $m_j=p_j-n_j'$. By Lemma~\ref{lem1} we have
\begin{align*}
J_j &\cong M_{p_j}(E_j)_{h_j} \cong M_{p_j}(E_j)_{d_j} \cong \\
&\cong
\left[\begin{array}{ccc} 0 & (1-d_j)M_{p_j}(E_j)d_j & (1-d_j)M_{p_j}(E_j)(1-d_j) \\  0 & d_jM_{p_j}(E_j)d_j & d_jM_{p_j}(E_j)(1-d_j) \\ 0 & 0 & 0 \\ \end{array}\right] \cong \\
&\cong
\left[\begin{array}{ccc} 0 & M_{m_j,n_j'}(E_j) & M_{m_j}(E_j) \\  0 & M_{n_j'}(E_j) & M_{n_j',m_j}(E_j) \\ 0 & 0 & 0 \\ \end{array}\right],
\end{align*}
where the last isomorphism is obvious due to the diagonal form of $d_j$. The Jacobson radical of $J_j$ is its strictly upper triangular part, so $J_j/J(J_j) \cong M_{n_j'}(E_j)$. The proof that $J_j$ is directly irreducible is the same as in the proof of Theorem~\ref{main_corner}, because $J_j$ is clearly generated by
$\left[\begin{array}{ccc} 0 & 0 & 0 \\  0 & 1 & 0 \\ 0 & 0 & 0 \\ \end{array}\right]$.
Since $aRa \cong J_0 \oplus \bigoplus_{j=1}^{k'} J_j$, the uniqueness in Theorem~\ref{main_corner} implies $k'=k$ and after a possible renumbering $J_j=I_j$, which additionally implies $n_j'=n_j$ and $E_j \cong D_j$ for all $1\leq j \leq k'=k$. For $j>k'=k$ we set $m_j=p_j$ and $D_j=E_j$.

Finally, observe that the right rank of $a$ is the same as the right rank of $e$, which is just $p_1+p_2\ldots+p_l=n_1+n_2+\ldots+n_k+m_1+m_2+\ldots+m_l$. \endproof

Note that the conclusion of Corollary~\ref{finite_corner} holds also for rings. We only needed the assumption that $R$ is an algebra to define $I_0$ as a direct complement of some ideal. If $a$ has finite rank this complement automatically exists. 

\begin{remark}
Let $a$ be an element of finite right rank in a semiprime algebra $R$. While the structure of $aRa$ uniquely determines the rank of $a^2$, it does not determine the rank of $a$. For example, take
\[a=\left[\begin{array}{cc} 0 & 1 \\ 0 & 0 \\ \end{array}\right] \in \left[\begin{array}{cc} \mathbb{H} & \mathbb{H} \\ \mathbb{H} & \mathbb{H} \\ \end{array}\right]=R
\quad \textup{and}\quad
b=\left[\begin{array}{cc} 0 & I \\ 0 & 0 \\ \end{array}\right] \in \left[\begin{array}{cc} M_2(\mathbb{R}) & M_2(\mathbb{R}) \\ M_2(\mathbb{R}) & M_2(\mathbb{R}) \\ \end{array}\right]=T,\]
where $I$ is the identity matrix. Then $aRa \cong bTb$ since these are both $4$-dimensional $\mathbb{R}$-algebras with trivial multiplication, however, $\rank a=1$ while $\rank b=2$. 
\end{remark}

Finally, we present an example, which shows that the assumption that $a^2$ has finite rank is crucial for the direct sum decomposition of $aRa$. There are two things that can go wrong if $a^2$ has infinite rank; the factor ring $aRa/J(aRa)$ need not be artinian and the potential decomposition of $aRa$ need not be direct. Even assuming that $a^3$ has finite rank is not sufficient.

\begin{example}
Let $E=\End_F(V)$, where $V$ is an infinite dimensional vector space over some field $F$, and let $I$ be the identity endomorphism of $V$. Consider the following element $a$ and algebra $R \subseteq M_{10}(E)$
{% for \arraycolsep
\arraycolsep=2.5pt
\[a=\left[
\begin{array}{cccccccccc}
 0 & I & 0 & 0 & 0 & 0 & 0 & 0 & 0 & 0 \\
 0 & 0 & 0 & 0 & 0 & 0 & 0 & 0 & 0 & 0 \\
 0 & 0 & 0 & I & 0 & 0 & 0 & 0 & 0 & 0 \\
 0 & 0 & 0 & 0 & I & 0 & 0 & 0 & 0 & 0 \\
 0 & 0 & 0 & 0 & 0 & 0 & 0 & 0 & 0 & 0 \\
 0 & 0 & 0 & 0 & 0 & 0 & I & 0 & 0 & 0 \\
 0 & 0 & 0 & 0 & 0 & 0 & 0 & I & 0 & 0 \\
 0 & 0 & 0 & 0 & 0 & 0 & 0 & 0 & 0 & 0 \\
 0 & 0 & 0 & 0 & 0 & 0 & 0 & 0 & 0 & I \\
 0 & 0 & 0 & 0 & 0 & 0 & 0 & 0 & 0 & 0
\end{array}
\right], \quad
R=\left[
\begin{array}{cccccccccc}
 E & E & E & E & E & E & E & E & E & E \\
 E & E & E & E & E & E & E & E & E & E \\
 0 & 0 & E & E & E & E & E & E & E & E \\
 0 & 0 & E & E & E & E & E & E & E & E \\
 0 & 0 & E & E & E & E & E & E & E & E \\
 0 & 0 & 0 & 0 & 0 & E & E & E & E & E \\
 0 & 0 & 0 & 0 & 0 & E & E & E & E & E \\
 0 & 0 & 0 & 0 & 0 & E & E & E & E & E \\
 0 & 0 & 0 & 0 & 0 & 0 & 0 & 0 & E & E \\
 0 & 0 & 0 & 0 & 0 & 0 & 0 & 0 & E & E
\end{array}
\right].\]
}% for \arraycolsep
Denote by $a^T$ the formal transpose of $a$ as a $10 \times 10$ matrix. Observe that $aa^Ta=a$, $a^2(a^2)^Ta^2=a^2$, and $a^3=0$, so that all powers of $a$ are regular and $a^3$ has finite rank. Clearly we have
{% for \arraycolsep
\arraycolsep=2.5pt
\[aRa=\left[
\begin{array}{cccccccccc}
 0 & E & 0 & E & E & 0 & E & E & 0 & E \\
 0 & 0 & 0 & 0 & 0 & 0 & 0 & 0 & 0 & 0 \\
 0 & 0 & 0 & E & E & 0 & E & E & 0 & E \\
 0 & 0 & 0 & E & E & 0 & E & E & 0 & E \\
 0 & 0 & 0 & 0 & 0 & 0 & 0 & 0 & 0 & 0 \\
 0 & 0 & 0 & 0 & 0 & 0 & E & E & 0 & E \\
 0 & 0 & 0 & 0 & 0 & 0 & E & E & 0 & E \\
 0 & 0 & 0 & 0 & 0 & 0 & 0 & 0 & 0 & 0 \\
 0 & 0 & 0 & 0 & 0 & 0 & 0 & 0 & 0 & E \\
 0 & 0 & 0 & 0 & 0 & 0 & 0 & 0 & 0 & 0
\end{array}
\right] \qquad\mathrm{and}\qquad aRa/J(aRa) \cong E \times E,\]
}% for \arraycolsep
so the factor algebra is not artinian. Nevertheless, following the proof of Theorem~\ref{main_corner}, we can extract a decomposition of $aRa$. Assuming the notation from Theorem~\ref{main_corner}, somewhat tedious computations (taking into account the isomorphism $aRa \cong (eRe)_f$) show, that we have $aRa=I_0+I_1+I_2$, where $I_0$, $I_1$, $I_2$ are consecutively equal to
{% for \arraycolsep
\arraycolsep=2.5pt
\begin{align*}
&\left[
\begin{array}{cccccccccc}
 0 & E & 0 & 0 & 0 & 0 & 0 & 0 & 0 & 0 \\
 0 & 0 & 0 & 0 & 0 & 0 & 0 & 0 & 0 & 0 \\
 0 & 0 & 0 & 0 & 0 & 0 & 0 & 0 & 0 & 0 \\
 0 & 0 & 0 & 0 & 0 & 0 & 0 & 0 & 0 & 0 \\
 0 & 0 & 0 & 0 & 0 & 0 & 0 & 0 & 0 & 0 \\
 0 & 0 & 0 & 0 & 0 & 0 & 0 & 0 & 0 & 0 \\
 0 & 0 & 0 & 0 & 0 & 0 & 0 & 0 & 0 & 0 \\
 0 & 0 & 0 & 0 & 0 & 0 & 0 & 0 & 0 & 0 \\
 0 & 0 & 0 & 0 & 0 & 0 & 0 & 0 & 0 & E \\
 0 & 0 & 0 & 0 & 0 & 0 & 0 & 0 & 0 & 0
\end{array}
\right],
\left[
\begin{array}{cccccccccc}
 0 & 0 & 0 & E & E & 0 & E & E & 0 & E \\
 0 & 0 & 0 & 0 & 0 & 0 & 0 & 0 & 0 & 0 \\
 0 & 0 & 0 & E & E & 0 & E & E & 0 & E \\
 0 & 0 & 0 & E & E & 0 & E & E & 0 & E \\
 0 & 0 & 0 & 0 & 0 & 0 & 0 & 0 & 0 & 0 \\
 0 & 0 & 0 & 0 & 0 & 0 & 0 & 0 & 0 & 0 \\
 0 & 0 & 0 & 0 & 0 & 0 & 0 & 0 & 0 & 0 \\
 0 & 0 & 0 & 0 & 0 & 0 & 0 & 0 & 0 & 0 \\
 0 & 0 & 0 & 0 & 0 & 0 & 0 & 0 & 0 & 0 \\
 0 & 0 & 0 & 0 & 0 & 0 & 0 & 0 & 0 & 0
\end{array}
\right],
\left[
\begin{array}{cccccccccc}
 0 & 0 & 0 & 0 & 0 & 0 & E & E & 0 & E \\
 0 & 0 & 0 & 0 & 0 & 0 & 0 & 0 & 0 & 0 \\
 0 & 0 & 0 & 0 & 0 & 0 & E & E & 0 & E \\
 0 & 0 & 0 & 0 & 0 & 0 & E & E & 0 & E \\
 0 & 0 & 0 & 0 & 0 & 0 & 0 & 0 & 0 & 0 \\
 0 & 0 & 0 & 0 & 0 & 0 & E & E & 0 & E \\
 0 & 0 & 0 & 0 & 0 & 0 & E & E & 0 & E \\
 0 & 0 & 0 & 0 & 0 & 0 & 0 & 0 & 0 & 0 \\
 0 & 0 & 0 & 0 & 0 & 0 & 0 & 0 & 0 & 0 \\
 0 & 0 & 0 & 0 & 0 & 0 & 0 & 0 & 0 & 0
\end{array}
\right].
\end{align*}
}% for \arraycolsep
However, the above sum is not direct because $I_1$ and $I_2$ intersect.
\end{example}


\begin{thebibliography}{99}
\bibitem{Bre-Sem} M. Bre\v sar, P. \v Semrl: Finite rank elements in semisimple Banach algebras, \emph{Studia Math.} 128 (1998), no. 3, 287--298.
\bibitem{Eis} D. Eisenbud: \emph{Commutative Algebra with a View Toward Algebraic Geometry}, Springer-Verlag, New York, 1995.
\bibitem{Lam} T.Y. Lam: \emph{A First Course in Noncommutative Rings, Second Edition}, Springer Science + Business Media, New York, 2001.
\bibitem{Lam2} T.Y. Lam: \emph{Lectures on Modules and Rings}, Springer-Verlag, New York, 1999.
\bibitem{Lam3} T.Y. Lam: Corner Ring Theory: A Generalization of Peirce Decompositions, I, \emph{Algebras, Rings and Their Representations: Proceedings of the International Conference on Algebras, Modules and Rings, Lisbon 2003}, World Scientific, 2006, pp. 153-182.
\bibitem{Sto1} N. Stopar: Rank of elements of general rings in connection with unit-regularity, preprint, arXiv:1809.06105v1 [math.RA].
\end{thebibliography}
\end{document}